\documentclass[a4paper,12pt]{amsart}
\usepackage{amsfonts}
\usepackage{amssymb}
\usepackage{ifthen}
\usepackage{graphicx}
\usepackage{color}
\nonstopmode \numberwithin{equation}{section}
\newtheorem{thm}{Theorem}
\newtheorem{lem}{Lemma}
\newtheorem{cor}{Corollary}

\newtheorem{cl}{Claim}
\newtheorem{ca}{Case}
\newtheorem{sca}{Subcase}
\newtheorem{scl}{Subclaim}
\newtheorem{conj}{Conjecture}

\theoremstyle{definition}
\newtheorem{defn}{Definition}

\newtheorem{op}[equation]{Open Problem}
\newtheorem{ques}[equation]{Question}
\newtheorem{rem}{Remark}[section]
\newtheorem{exam}[equation]{Example}

\newcounter {own}
\def\theown {\thesection       .\arabic{own}}

\newenvironment{pf}[1][]{%
 \vskip 3mm
 \noindent
 \ifthenelse{\equal{#1}{}}%
  {{\slshape Proof. }}%
  {{\slshape #1.} }%
 }%
{\qed\bigskip}

\newcounter{alphabet}
\newcounter{tmp}

\makeatletter
\newcommand{\Ref}[1]{\@ifundefined{r@#1}{}{\setcounter{tmp}{\ref{#1}}\Alph{tmp}}}
\makeatother

\newenvironment{Lem}[1][]{\refstepcounter{alphabet}%
\bigskip%
\noindent%
{\bf Lemma \Alph{alphabet}}%
{\bf .} \itshape}{\vskip 8pt}

\newcommand{\IC}{{\mathbb C}}
\newcommand{\ID}{{\mathbb D}}
\newcommand{\IB}{{\mathbb B}}

\def\be{\begin{equation}}
\def\ee{\end{equation}}

\newcommand{\bee}{\begin{enumerate}}
\newcommand{\eee}{\end{enumerate}}

\newcommand{\blem}{\begin{lem}}
\newcommand{\elem}{\end{lem}}
\newcommand{\bthm}{\begin{thm}}
\newcommand{\ethm}{\end{thm}}
\newcommand{\bcor}{\begin{cor}}
\newcommand{\ecor}{\end{cor}}
\newcommand{\beg}{\begin{exam}}
\newcommand{\eeg}{\end{exam}}
\newcommand{\begs}{\begin{examples}}
\newcommand{\eegs}{\end{examples}}
\newcommand{\bdefe}{\begin{defn}}
\newcommand{\edefe}{\end{defn}}
\newcommand{\bprob}{\begin{prob}}
\newcommand{\eprob}{\end{prob}}
\newcommand{\bques}{\begin{ques}}
\newcommand{\eques}{\end{ques}}
\newcommand{\bei}{\begin{itemize}}
\newcommand{\eei}{\end{itemize}}
\newcommand{\bcon}{\begin{conj}}
\newcommand{\econ}{\end{conj}}
\newcommand{\bop}{\begin{op}}
\newcommand{\eop}{\end{op}}

\newcommand{\bca}{\begin{ca}}
\newcommand{\eca}{\end{ca}}
\newcommand{\bsca}{\begin{sca}}
\newcommand{\esca}{\end{sca}}

\newcommand{\bcl}{\begin{cl}}
\newcommand{\ecl}{\end{cl}}

\newcommand{\bscl}{\begin{scl}}
\newcommand{\escl}{\end{scl}}

\newcommand{\bcons}{\begin{conjs}}
\newcommand{\econs}{\end{conjs}}
\newcommand{\bprop}{\begin{propo}}
\newcommand{\eprop}{\end{propo}}
\newcommand{\br}{\begin{rem}}
\newcommand{\er}{\end{rem}}
\newcommand{\brs}{\begin{rems}}
\newcommand{\ers}{\end{rems}}
\newcommand{\bo}{\begin{obser}}
\newcommand{\eo}{\end{obser}}
\newcommand{\bos}{\begin{obsers}}
\newcommand{\eos}{\end{obsers}}
\newcommand{\bpf}{\begin{pf}}
\newcommand{\epf}{\end{pf}}
\newcommand{\ba}{\begin{array}}
\newcommand{\ea}{\end{array}}
\newcommand{\beq}{\begin{eqnarray}}
\newcommand{\beqq}{\begin{eqnarray*}}
\newcommand{\eeq}{\end{eqnarray}}
\newcommand{\eeqq}{\end{eqnarray*}}

\newcommand{\ra}{\rightarrow}

\newcommand{\ds}{\displaystyle}

\newcounter{minutes}
\divide\time by 60
\newcounter{hours}
\multiply\time by 60 \addtocounter{minutes}{-\time}

\begin{document}

\bibliographystyle{amsplain}
\title{Distortion and covering theorems  of pluriharmonic mappings
}

\def\thefootnote{}
\footnotetext{ \texttt{\tiny File:~\jobname .tex,
          printed: \number\day-\number\month-\number\year,
          \thehours.\ifnum\theminutes<10{0}\fi\theminutes}
} \makeatletter\def\thefootnote{\@arabic\c@footnote}\makeatother

\author{Sh. Chen }
\address{Sh. Chen, Department of Mathematics and Computational
Science, Hengyang Normal University, Hengyang, Hunan 421008,
People's Republic of China.} \email{mathechen@126.com}

\author{S. Ponnusamy $^\dagger $}
\address{S. Ponnusamy,
Indian Statistical Institute (ISI), Chennai Centre, SETS (Society
for Electronic Transactions and security), MGR Knowledge City, CIT
Campus, Taramani, Chennai 600 113, India. }
\email{samy@isichennai.res.in, samy@iitm.ac.in}


\subjclass[2000]{Primary: 30C99, 31C10; Secondary: 15A04, 30C62, 32A10.}
\keywords{Pluriharmonic  mapping, distortion theorem, linear invariant family, quasiregular mapping.
\\ $^\dagger$ {Corresponding author.}
}



\begin{abstract}
In this paper, we mainly investigate distortion and covering
theorems on some classes of pluriharmonic mappings.
\end{abstract}


\maketitle \pagestyle{myheadings} \markboth{Pluriharmonic mappings}{Sh. Chen and S. Ponnusamy}

\section{Introduction and preliminaries}\label{csw-sec1}

The notion of linear-invariant family (hereafter $\mathcal{LIF}$) of holomorphic functions defined on the unit disk $
\ID := \{z\in \IC:\, |z| <1 \}$ was first introduced  by Pommerenke in  \cite{P} and showed
a number of important properties of such families.
Recall that if $\mathcal A$ denotes the family of all holomorphic functions $f$
on $\ID$ with the topology of uniform convergence of compact subsets of $\ID$, then a
subfamily $\mathcal F$ of $\mathcal A$ is called linear-invariant if it is closed
under the re-normalized composition with a conformal automorphism of $\ID$.
If the modulus of the second Taylor coefficient is bounded in $\mathcal F$, then the order $\alpha$
of the $\mathcal{LIF}$ is defined to be
$$\alpha := \sup\{|f''(0)|/2:\, f\in {\mathcal F}\}.
$$
Many properties of a $\mathcal{LIF}$ depends on the order of the
family. A universal $\mathcal{LIF}$ of order $\alpha$, denoted by
${\mathcal U}_\alpha$, is the union of all $\mathcal{LIF}$'s
$\mathcal F$ such that order of $\mathcal F$ less than or equal to
$\alpha$. The fact is that ${\mathcal U}_\alpha$ is empty if $\alpha
<1$ and ${\mathcal U}_1$ coincides with the  family of all
normalized holomorphic functions $f$ which univalently map $\ID$
onto convex domains, see \cite{P}. Also, a $\mathcal{LIF}$ of order
$2$ is the family $\mathcal S$ of normalized univalent functions
from $\mathcal A$. Moreover, it  has been proved that many
subfamilies of univalent mappings on $\ID$ are linearly invariant,
see for example \cite{koepf} and the references therein. For the
regularity growth of functions on ${\mathcal U}_\alpha$, we refer to
\cite{Cam,S1a,S1b}. The concept of linear invariance was generalized
by many authors in many different contexts and in 1997, Pfaltzgraff
\cite{Pf} extended  this concept for locally holomorphic functions
defined on the unit ball of the complex Euclidean $n$-space
$\mathbb{C}^{n}$ and many properties were further discussed in
\cite{PS,PS2}. For our discussion, we need to deal with such
problems in the higher dimensional case.

As with the standard practice, for $z=(z_{1}~\cdots~z_{n})$ and
$w=(w_{1}~\cdots~w_{n})$ in $\mathbb{C}^{n}$, we let
$\overline{z}=(\overline{z}_{1}~\cdots~\overline{z}_{n} ),$ and
$\langle z,w\rangle := \sum_{k=1}^nz_k\overline{w}_k$ with the
associated Euclidean norm $ \|z\|:={\langle z,z\rangle}^{1/2}$ which
makes $\mathbb{C}^n$ into an $n$-dimensional complex Hilbert space.
Throughout the discussion an element  $z\in \mathbb{C}^{n}$  is identified as an $n\times 1$
column vector.
For $a\in \mathbb{C}^n$ and $r>0$,
$$\IB^n(a,r)=\{z\in \mathbb{C}^{n}:\, \|z-a\|<r\}
$$
denotes the (open) ball of radius $r$ with center $a$. Also, we let
$\IB^n(r):=\IB^n(0,r)$ and use $\IB^n$ to denote the unit ball
$\IB^n(1)$, and $\mathbb{D}=\mathbb{B}^1$.

A continuous complex-valued function $f$ defined on a domain
$G\subset\mathbb{C}^{n}$ is said to be {\it pluriharmonic} if for
each fixed $z\in G$ and $\theta\in\partial\mathbb{B}^{n}$, the
function $f(z+\theta\zeta)$ is harmonic in $\{\zeta\in \IC:\;
\|\theta\zeta-z\|< d_{G}(z)\}$, where $d_{G}(z)$ denotes the
distance from $z$ to the boundary $\partial G$ of $G$. It follows
from \cite[Theorem 4.4.9]{R} that a real-valued function $u$ defined
on $G$ is pluriharmonic if and only if it is locally the real part
of a holomorphic function. If $\Omega$ is a simply connected domain
in $\mathbb{C}^{n}$, then it is clear that a mapping $f:\,\Omega\ra
\mathbb{C}$ is pluriharmonic if and only if $f$ has a representation
$f=h+\overline{g}$, where  $h, g$ are holomorphic in $\Omega$ (cf.
\cite{Vl}). A {\it vector-valued mapping}
$f=(f_{1}~\cdots~f_{N})^{T}$
defined in $\mathbb{B}^{n}$ is said to be  pluriharmonic, if each
component $f_{j}$ ($1\leq j\leq N$) is a pluriharmonic mapping from
$\mathbb{B}^{n}$ into $\mathbb{C}$, where $N$ is a positive integer
and $T$ is the transpose of a matrix.
 We refer to
\cite{CPW-1,CPW-2,CPW-3,DHK,I,R} for further details and recent
investigations on pluriharmonic mappings.

For an $n\times n$ complex matrix $A$, we
introduce the {\it operator norm}
$$\|A\|=\sup_{z\neq0}\frac{\|Az\|}{\|z\|}=\max\left\{\|A\theta\|:\,\theta\in\partial\mathbb{B}^{n}\right\}.
$$

We use $L(\mathbb{C}^{n},\mathbb{C}^{m})$ to denote the space of
continuous {\it linear operators} from $\mathbb{C}^{n}$ into
$\mathbb{C}^{m}$ with the operator norm, and let $I_{n}$ be the {\it
identity operator} in $L(\mathbb{C}^{n},\mathbb{C}^{n})$.

We denote by $\mathcal{PH}(\mathbb{B}^{n}, \mathbb{C}^n)$ the set of
all {\it vector-valued pluriharmonic mappings} from $\mathbb{B}^{n}$
into $\mathbb{C}^n$. Then every $f\in\mathcal{PH}(\mathbb{B}^{n}, \mathbb{C}^n)$ can be written as
$f=h+\overline{g}$, where $h$ and $g$ are holomorphic in $\mathbb{B}^{n}$,
and this representation is unique when $g(0) = 0$. It is a simple exercise to see that the real
Jacobian  determinant of $f$ can be written as
$$\det J_{f}=\det\left(\begin{array}{cccc}
\ds Dh & \overline{Dg}  \\
\ds Dg & \overline{Dh}
\end{array}\right)
$$
and  if $h$ is locally biholomorphic (i.e. the complex Jacobian
matrix $J_{f}(z)$  of $f$ at each $z$ is invertible), then the
determinant of $J_{f}$ has the form \be\label{eq-ex2a} \det
J_{f}=|\det
Dh|^{2}\det\left(I_{n}-Dg[Dh]^{-1}\overline{Dg[Dh]^{-1}}\right).\ee
In the case
of a {\it planar harmonic mapping} $f=h+\overline{g}$, we find that
$$\det J_f=|h'|^2-|g'|^2,
$$
and so, $f$ is locally univalent and sense-preserving in $\ID$ if
and only if $|g'(z)|<|h'(z)|$ in $\ID$; or equivalently if $h'(z)
\neq0$ and the dilatation $\omega (z)=g'(z)/h'(z)$ is analytic in
$\ID$ and has the property that $|\omega (z)|<1$ in $\ID$ (see
\cite{Du,Lewy}). For
$f=h+\overline{g}\in\mathcal{PH}(\mathbb{B}^{n}, \mathbb{C}^n)$, the
condition $\|Dg[Dh]^{-1}\|<1$  is sufficient for $\det J_f$ to be
positive and hence for $f$ to be sense-preserving. This is indeed a
natural generalization of one-variable condition (cf. \cite{DHK}).

For motivation, consider the Taylor expansion of a function
$f=h+\overline{g}\in\mathcal{PH}(\mathbb{B}^{n}, \mathbb{C}^n)$ with
$h(0)=g(0)=0$, where
\beq\label{eqt-1}
h(z)&=&[Dh(0)]z+\frac{1}{2}[D^{2}h(0)](z,z)+\cdots+\frac{1}{m}[D^{m}h(0)](z,\ldots,z)+\cdots\\
\nonumber &=&A_{1}z+A_{2}(z,z)+A_{m}(z,\ldots,z)+\cdots
\eeq
and
\beq\label{eqt-2}
g(z)&=&[Dg(0)]z+\frac{1}{2}[D^{2}g(0)](z,z)+\cdots+\frac{1}{m}[D^{m}g(0)](z,\ldots,z)+\cdots\\
\nonumber &=&B_{1}z+B_{2}(z,z)+B_{m}(z,\ldots,z)+\cdots.
\eeq

As with one variable case, a $\mathcal{LIF}$ in
$\mathbb{B}^{n}$ is a family $\mathcal{M}$ of locally biholomorphic
mappings $f:\,\mathbb{B}^{n}\rightarrow\mathbb{C}^{n}$ such that if
$f\in\mathcal{M}$ then
\begin{enumerate}
\item[(i)] $f(0)=0,$ $Df(0)=I_{n}$ and
\item[(ii)] $\Lambda_{\phi}(f)\in\mathcal{M}$ for all
$\phi\in\mbox{Aut}(\mathbb{B}^{n})$, the holomorphic automorphism of
$\mathbb{B}^{n}$.
\end{enumerate}
Here $\Lambda_{\phi}(f)=[D\phi(0)]^{-1}[Df(\phi(0))]^{-1}[f(\phi(z))-f(\phi(0))]$
denotes the {\it Koebe transform} of $f$ (cf. \cite{PS,PS2}) and thus, the classical definition of the order
$\alpha$ of $\mathcal{LIF}$ introduced in the beginning is generalized as follows:

\begin{defn}
If $\mathcal{M}$ is a $\mathcal{LIF}$, then the {\it norm order} of
$\mathcal{M}$ is the quantity
$$\|{\rm ord}\|_{\mathcal{M}}=\sup\left\{\frac{1}{2}\|D^{2}f(0)\|:\,f\in\mathcal{M}\right\}=\alpha.
$$
\end{defn}

In \cite[Theorem 3.1]{PS}, it has been shown that $\alpha\geq1$. As in the planar case, the universal
linearly-invariant family $\mathcal{M}_{\alpha}$ of order $\alpha$
is defined as the union of all linearly invariant families of order
less than or equal to $\alpha$ (cf. \cite{P}).

Our main aim of this paper is to extend the corresponding results of
 \cite{S1} and \cite{S2} to higher dimensional case.

\section{Main results}


Let $\mathcal{PH}(\alpha,k)$ denote the set of all sense-preserving
mappings $f=h+\overline{g}\in\mathcal{PH}(\mathbb{B}^{n},
\mathbb{C}^n)$ with the normalization $h(0)=g(0)=0$,
$\|Dh(0)+\overline{Dg(0)}\|=1,$ $[Dh(0)]^{-1}h(z)\in
\mathcal{M}_{\alpha}$, and such that for $k\in[0,1)$,
$$\left\|Dg(z)[Dh(z)]^{-1}\right\|\leq k,
$$
where $h$ is locally biholomorphic  and $g$ is holomorphic in $\mathbb{B}^{n}$.


Obviously, if $n=1$, then $\mathcal{PH}(\alpha,k)$ coincides with
the set $H(\alpha, K)$ of \cite{S1} and \cite{S2}. As a
generalization of \cite[Theorem 1]{S1}, we have.

\begin{thm}\label{thm-1}
For $\alpha<\infty$, the classes $\mathcal{PH}(\alpha,k)$ are
compact with respect to the topology of almost uniform convergence
in $\mathbb{B}^{n}$.
\end{thm}


The derivative of $f=h+\overline{g}\in\mathcal{PH}(\mathbb{B}^{n},
\mathbb{C}^n)$ in the direction of vector
$\theta\in\partial\mathbb{B}^{n}$ at the point $z$ will be denoted
by
$$\partial_{\theta}f(z)=\lim_{\rho\rightarrow0+}\frac{f(z+\rho\theta)-f(z)}{\rho}=Dh(z)\theta+\overline{Dg(z)\theta},
$$
where $h$ and $g$ are holomorphic in $\mathbb{B}^{n}.$ We use the standard notations:
$$ \Lambda_{f}=\max_{\theta\in\partial\mathbb{B}^{n}}\left\|\partial_{\theta}f\right\|\;\;
\mbox{ and }\;\; \lambda_{f}=\min_{ \theta\in\partial\mathbb{B}^{n}}\left\|\partial_{\theta}f\right\|.
$$
With this setting, we now present a generalization of \cite[Theorem 2]{S1}.

\begin{thm}\label{thm-2}
For $\alpha<\infty$, let
$f=h+\overline{g}\in\mathcal{PH}(\alpha,k)$. Then
\be\label{eq-2r}\frac{1-k}{\left\|[Dh(0)]^{-1}\right\|}\frac{(1-\|z\|)^{\alpha-1}}{(1+\|z\|)^{\alpha+1}}\leq\Lambda_{f}(z)
\leq \left (\frac{1+k}{1-k}\right
)\frac{(1+\|z\|)^{\alpha-1}}{(1-\|z\|)^{\alpha+1}} \ee and
\be\label{eq-2.1r}\|f(z)\|\leq\frac{1+k}{2\alpha(1-k)}\left\{\frac{(1+\|z\|)^{\alpha}}{(1-\|z\|)^{\alpha}}-1\right\}.
\ee
 In particular,
if $n=1$, then the estimate of {\rm(\ref{eq-2r})} is sharp for
$\theta=\pm\frac{\pi}{2}$. Moreover, if $z=re^{it}$,  then the
equality on the right of {\rm(\ref{eq-2r})} is obtained for
$f(z)=h(z)-k\overline{h(z)},$
$$h(z)=\frac{e^{it}}{2\alpha(1-k)}\left[\left(\frac{1+ze^{-it}}{1-ze^{-it}}\right )^{\alpha}-1\right]
$$
and the equality on the left of {\rm(\ref{eq-2r})} is obtained for $f(z)=h^{\ast}(z)+k\overline{h^{\ast}(z)}$,
$$h^{\ast}(z)=\frac{e^{it}}{2\alpha(1+k)}\left[\left (\frac{1-ze^{-it}}{1+ze^{-it}}\right )^{\alpha}-1\right].
$$
\end{thm}



The following result is a covering theorem of
$\mathcal{PH}(\alpha,k)$.

\begin{thm}\label{thm-4}
For $r\in(0,1]$ and $\alpha<\infty$, if
$f=h+\overline{g}\in\mathcal{PH}(\alpha,k)$, then
$f(\mathbb{B}^{n}(r))$ contains a univalent ball $\mathbb{B}^{n}(R)$
with
$$R\geq\frac{(1-k)|\det Dh(0)|}{\|Dh(0)\|^{n-1}}\int_{0}^{r}
\frac{(1-x)^{(2n-1)\alpha+(n-3)/2}}{(1+x)^{(2n-1)\alpha-(n-3)/2}}\,dx.
$$
In particular, if $n=1$, then
$R=(1-k)\left[1-\big(\frac{1-r}{1+r}\big)^{\alpha}\right]/[2\alpha(1+k)],$
and  the extreme function  $f=h+k\overline{h}$ shows that this
estimate is sharp, where
$$h(z)=\frac{\pm i}{2\alpha(1+k)}\left[\left(\frac{1\pm iz}{1\mp iz}\right )-1\right].
$$
\end{thm}

We remark that Theorem \ref{thm-4} is a generalization of \cite[Theorem 3]{S1}.

\begin{thm}\label{thm-5}
For $\alpha<\infty$, if
$f=h+\overline{g}\in\mathcal{PH}(\alpha,k)$, then
$$|\det J_{f}(z)|\geq \frac{(1-k^{2})^{n}}{\big(\det [Dh(0)]^{-1}\big)^{2}}\frac{\big(1-\|z\|\big)^{2n\alpha-n-1}}
{\big(1+\|z\|\big)^{2n\alpha+n+1}}.
$$
\end{thm}



For $r\in(0,1)$, a univalent mapping $f=h+\overline{g}\in\mathcal{PH}(\mathbb{B}^{n}, \mathbb{C}^n)$ with
$h(0)=g(0)=0$, $Dg(0)=0$ and
$$\big\|Dg[Dh]^{-1}\big\|<1
$$
is called {\it fully starlike} if it maps every ball
$\overline{\mathbb{B}^{n}(r)}$ onto a starlike domain with respect
to the origin, where $h$ is locally biholomorphic and $g$ is
holomorphic in $\mathbb{B}^{n}$ (cf. \cite{CDO}). The following
result is a generalization of \cite[Theorem 1.3]{CPW-4}.

\begin{thm}\label{th-7}
Let $r\in(0,1)$ and $f=h+\overline{g}\in\mathcal{PH}(\mathbb{B}^{n}, \mathbb{C}^n)$ be fully starlike, where $h$ is
locally biholomorphic and $g$ is holomorphic in $\mathbb{B}^{n}$. Then for all $z\in\overline{\mathbb{B}^{n}(r)}$,
$$\|h(z)\|\leq\frac{1}{1-r}\|f(z)\|.
$$
Furthermore, if $h\in\mathcal{M}_{\alpha}$, then
\begin{enumerate}
\item[\rm(a)] for $z\in\mathbb{B}^{n}(r_{0})$,
$$\|f(z)\|\geq r_{0}^{2}(1-r_{0})\frac{\|z\|}{(r_{0}+\|z\|)^{2}},
$$
where $r_{0}=4\alpha/(1+4\alpha^{2});$

\item[{\rm(b)}] $f$ differs from zero in $\mathbb{B}^{n}(r_{0})\backslash\{0\}$.
\end{enumerate}
\end{thm}



We remark that
$$\frac{4\alpha}{1+4\alpha^{2}}=\frac{1}{\alpha}-\frac{1}{\alpha(1+4\alpha^{2})}\thicksim\frac{1}{\alpha}
$$
as $\alpha\rightarrow\infty.$ Hence Theorem \ref{th-7}(b) is a generalization of \cite[Theorem 1]{S2}.

A continuous mapping $f:\ \Omega\subset\mathbb{R}^{n}\rightarrow\mathbb{R}^{n}$ is called {\it
$K$-quasiregular} if $f\in W_{n,{\rm loc}}^{1}(\Omega)$ and
$$\|Df(x)\|^{n}\leq K\det J_{f}(x)\ \mbox{for almost every}\ x\in\Omega,
$$
where $K$ $(\geq1)$ is a constant. Here $f\in W_{n,{\rm loc}}^{1}(\Omega)$ means that the distributional derivatives
$\partial f_{j}/\partial x_{k}$ of the coordinates $f_{j}$ of $f$ are locally in $L^{n}(\Omega)$ and $J_{f}(x)$
denotes the Jacobian of $f$ (cf. \cite{V}).

Let $f=(f_{1}~\cdots~ f_{n})^{T}\in \mathcal{PH}(\mathbb{B}^{n},
\mathbb{C}^n)$. For $j\in\{1,\ldots,n\}$, we let $z=(z_{1}~\cdots~
z_{n})^{T}$, $z_{j}=x_{j}+iy_{j}$ and $f_{j}(z)=u_{j}(z)+iv_{j}(z)$,
where $u_{j}$ and $v_{j}$ are real pluriharmonic functions from
$\mathbb{B}^{n}$ into $\mathbb{R}$. We denote the real Jacobian
matrix of $f$ by
$$J_{f}=\left(\begin{array}{cccc}
\ds \frac{\partial u_{1}}{\partial x_{1}}\; \frac{\partial
u_{1}}{\partial y_{1}}\; \frac{\partial u_{1}}{\partial x_{2}}\;
\frac{\partial u_{1}}{\partial y_{2}}\;\cdots\;
 \frac{\partial u_{1}}{\partial x_{n}}\; \frac{\partial u_{1}}{\partial y_{n}}\\[4mm]
 \ds \frac{\partial v_{1}}{\partial x_{1}}\; \frac{\partial v_{1}}{\partial y_{1}}\;
\frac{\partial v_{1}}{\partial x_{2}}\; \frac{\partial
v_{1}}{\partial y_{2}}\;\cdots\;
 \frac{\partial v_{1}}{\partial x_{n}}\; \frac{\partial v_{1}}{\partial  y_{n}}\\[2mm]
\vdots \\[2mm]
 \ds \frac{\partial u_{n}}{\partial x_{1}}\; \frac{\partial u_{n}}{\partial y_{1}}\;
\frac{\partial u_{n}}{\partial x_{2}}\; \frac{\partial
u_{n}}{\partial y_{2}}\;\cdots\;
 \frac{\partial u_{n}}{\partial x_{n}}\; \frac{\partial u_{n}}{\partial y_{n}}\\[4mm]
\ds  \frac{\partial v_{n}}{\partial x_{1}}\; \frac{\partial
v_{n}}{\partial y_{1}}\; \frac{\partial v_{n}}{\partial x_{2}}\;
\frac{\partial v_{n}}{\partial y_{2}}\;\cdots\;
\frac{\partial v_{n}}{\partial x_{n}}\; \frac{\partial v_{n}}{\partial y_{n}}
\end{array}\right).
$$
Let $\mathbb{B}^{2n}_{\mathbb{R}}$ denote the unit ball of
$\mathbb{R}^{2n}$. Then
$$\Lambda_{f}=\max_{\theta\in\partial\mathbb{B}^{2n}_{\mathbb{R}}}\|J_{f}\theta\|\
\mbox{and}\
\lambda_{f}=\min_{\theta\in\partial\mathbb{B}^{2n}_{\mathbb{R}}}\|J_{f}\theta\|.
$$

\begin{thm}\label{thm-6}
Let $f=h+\overline{g}\in \mathcal{PH}(\mathbb{B}^{n}, \mathbb{C}^n)$
with $\big\|Dg(z)[Dh(z)]^{-1}\big\|\leq c<1$ for
$z\in\mathbb{B}^{n}$, where $c$ is a positive constant. Then

\begin{enumerate}
\item[\rm(a)]  $f$ is a quasiregular mapping if and only if $h$ is a
quasiregular mapping;

\item[\rm(b)] $f(\mathbb{B}^{n})$ contains a univalent ball with the
radius
$$R\geq\frac{k_{n}\pi}{8m}\left(\frac{k_{n}\pi\sqrt{1-c}}{4K\sqrt{1+c}\log(1/(1-k_{n}))}\right)^{4n-1},
$$
where $m\approx4.2$, $\det J_{f}(0)=1$, $h$ is a $K$-quasiregular
mapping with $K\geq1$ and $0<k_{n}<1$ is a unique root such that
\be\label{eq-ex1} -4n\log(1-k_{n})=(4n-1)\frac{k_{n}}{1-k_{n}}. \ee
\end{enumerate}
\end{thm}

The roots $k_{n}$ in $(0,1)$ of the equation \eqref{eq-ex1} for the values of $n=1,2,3,4, 5$ are listed in Table \ref{tab1}
for a ready reference.
\begin{table}
\begin{tabular}{|l|l|}
  \hline
  Value of $n$ & Value of $k_n$ \\
  \hline
  1 & 0.423166 \\
  \hline
  2 & 0.230006 \\
  \hline
  3 & 0.157659 \\
  \hline
  4 & 0.119898 \\
  \hline
  5 & 0.0967215 \\
  \hline
\end{tabular}
\vspace{.2cm}
\caption{Values of $k_{n}$ of Equation \eqref{eq-ex1} for $n=1,2,3,4, 5$ \label{tab1}}
\end{table}


The proofs of Theorems \ref{thm-1}$-$\ref{thm-6} will be presented
in Section \ref{csw-sec2}.

\section{Proofs of the main theorems }\label{csw-sec2}
\subsection*{Proof of Theorem \ref{thm-1}}
 Consider a sequence $f_{m}=h_{m}+\overline{g}_{m}\in\mathcal{PH}(\alpha,k).$ By definition, we have the conditions
$\|Dh_{m}(0)+\overline{Dg_{m}(0)}\|=1$ and $\left\|Dg_{m}(z)[Dh_{m}(z)]^{-1}\right\|\leq k,$ we see that
$$\|Dh_{m}(0)\|\leq 1+\|Dg_{m}(0)\|
$$
whereas the second condition gives
$$ \|Dg_{m}(0)\|=\left\|Dg_{m}(0)[Dh_{m}(0)]^{-1}[Dh_{m}(0)]\right\|
\leq k\|Dh_{m}(0)\|.
$$
Using the last two inequalities, we easily have
\be\label{eq-1.1}
\|Dg_{m}(0)\|\leq\frac{k}{1-k}~\mbox{ and }~\|Dh_{m}(0)\|\leq\frac{1}{1-k}.
\ee
By (\ref{eq-1.1}), $[Dh_{m}(0)]^{-1}h_{m}(z)\in \mathcal{M}_{\alpha}$ and thus by \cite[Theorem 4.1]{PS}, we obtain that
\be\label{eq-1.2}
\frac{(1-\|z\|)^{\alpha-1}}{(1+\|z\|)^{\alpha+1}}\leq\left\|[Dh_{m}(0)]^{-1}Dh_{m}(z)\right\|\leq
\frac{(1+\|z\|)^{\alpha-1}}{(1-\|z\|)^{\alpha+1}},
\ee
which implies
\begin{eqnarray*}
\|[Dh_{m}(z)]\|&=&\left\|Dh_{m}(0)[Dh_{m}(0)]^{-1}Dh_{m}(z)\right\|\\
&\leq&\left\|[Dh_{m}(0)]^{-1}Dh_{m}(z)\right\|\|Dh_{m}(0)\|\\
&\leq&\frac{1}{(1-k)}\frac{(1+\|z\|)^{\alpha-1}}{(1-\|z\|)^{\alpha+1}}.
\end{eqnarray*}
Moreover, by the definition of $\mathcal{PH}(\alpha,k)$, it follows that
$$\|Dg_{m}(z)\|\leq k\|Dh_{m}(z)\|\leq\frac{k}{(1-k)}\frac{(1+\|z\|)^{\alpha-1}}{(1-\|z\|)^{\alpha+1}}.
$$
Hence $Dh_{m}(z)$ and $Dg_{m}(z)$ are uniformly bounded in compact
subsets of $\mathbb{B}^{n}$, which implies $\mathcal{PH}(\alpha,k)$
are compact. \hfill $\Box$

\subsection*{Proof of Theorem \ref{thm-2}} Let $f=h+\overline{g}\in \mathcal{PH}(\alpha,k)$ for some $\alpha<\infty$.
By the definition of directional derivatives, we have
\begin{eqnarray*}
\left\|\partial_{\theta}f(z)\right\|
&=&\left\|Dh(z)\theta+\overline{Dg(z)[Dh(z)]^{-1}Dh(z)\theta}\right\|\\
&\geq&\|Dh(z)\theta\|\left(1-\big\|Dg(z)[Dh(z)]^{-1}\big\|\right)\\
&\geq&(1-k)\|Dh(z)\theta\|
\end{eqnarray*}
and similarly,
\begin{eqnarray*}
\left\|\partial_{\theta}f(z)\right\|
&\leq&\|Dh(z)\theta\|\left(1+\big\|Dg(z)[Dh(z)]^{-1}\big\|\right)\\
&\leq&(1+k)\|Dh(z)\theta\|.
\end{eqnarray*}
It follows that
\be\label{eq-3}
(1-k)\|Dh(z)\|\leq\Lambda_{f}(z)=\max_{\theta\in\partial\mathbb{B}^{n}}
\left\|\partial_{\theta}f(z)\right\|\leq(1+k)\|Dh(z)\|.
\ee
Again, by elementary calculations, we have
$$\|Dh(z)\|=\left\|Dh(0)[Dh(0)]^{-1}Dh(z)\right\|\leq\left\|[Dh(0)]^{-1}Dh(z)\right\|\|Dh(0)\|,
$$
which gives
\be\label{eq-1.3}\frac{\|Dh(z)\|}{\|Dh(0)\|}\leq\left\|[Dh(0)]^{-1}Dh(z)\right\|\leq\|Dh(z)\|\left\|[Dh(0)]^{-1}\right\|.
\ee
By  $[Dh(0)]^{-1}h(z)\in \mathcal{M}_{\alpha}$ and \cite[Theorem 4.1]{PS}, we deduce that
\be\label{eq-1.4}
\frac{(1-\|z\|)^{\alpha-1}}{(1+\|z\|)^{\alpha+1}}\leq\left\|[Dh(0)]^{-1}Dh(z)\right\|\leq
\frac{(1+\|z\|)^{\alpha-1}}{(1-\|z\|)^{\alpha+1}}.
\ee
By (\ref{eq-1.3}) and (\ref{eq-1.4}), we get
\be\label{eq-1.5}\frac{1}{\left\|[Dh(0)]^{-1}\right\|}\frac{(1-\|z\|)^{\alpha-1}}{(1+\|z\|)^{\alpha+1}}\leq\|Dh(z)\|\leq
\frac{(1+\|z\|)^{\alpha-1}}{(1-\|z\|)^{\alpha+1}}\|Dh(0)\|,\ee which
implies
\be\label{eq-1.6}
\frac{1-k}{\left\|[Dh(0)]^{-1}\right\|}\frac{(1-\|z\|)^{\alpha-1}}{(1+\|z\|)^{\alpha+1}}\leq\Lambda_{f}(z)
\leq \frac{(1+\|z\|)^{\alpha-1}}{(1-\|z\|)^{\alpha+1}}\|Dh(0)\|(1+k).
\ee
Applying (\ref{eq-1.6}) and the inequality,
\be\label{eqt}
\frac{1}{1+k}\leq\|Dh(0)\|\leq\frac{1}{1-k},
\ee
we conclude that
\be\label{eq2.2r}\frac{1-k}{\left\|[Dh(0)]^{-1}\right\|}\frac{(1-\|z\|)^{\alpha-1}}{(1+\|z\|)^{\alpha+1}}\leq\Lambda_{f}(z)
\leq
\frac{1+k}{(1-k)}\frac{(1+\|z\|)^{\alpha-1}}{(1-\|z\|)^{\alpha+1}}.
\ee

Now we prove (\ref{eq-2.1r}). Let $[0,z]$ be the segment from $0$ to
$z\in\mathbb{B}^{n}$. Then by using (\ref{eq2.2r}), we have   

\begin{eqnarray*}
\|f(z)\|&=&\left\|\int_{[0,z]}df(\zeta)\right\|=\left\|\int_{[0,z]}Dh(\zeta)\,d\zeta+\overline{Dg(\zeta)\,d\zeta}\right\|\\
&\leq&\int_{[0,z]}\Lambda_{f}(\zeta)\|d\zeta\|\\
&=&\frac{1+k}{1-k}\int_{0}^{1}\frac{(1+t\|z\|)^{\alpha-1}}{(1-t\|z\|)^{\alpha+1}}\|z\|\,dt\\
&=&\frac{1+k}{2\alpha(1-k)}\left\{\frac{(1+\|z\|)^{\alpha}}{(1-\|z\|)^{\alpha}}-1\right\}.
\end{eqnarray*}
The proof of this theorem is complete. \hfill $\Box$

\begin{Lem} {\rm (\cite[Lemma 4]{LX})}\label{Lem-A}
Let $A$ be an $n \times n$ complex $($real$)$ matrix with $\|A\|\neq0$. Then for all unit vector
$\theta\in\partial \mathbb{B}^{n}$, the inequality
$$\|A\theta\|\geq\frac{|\det A|}{\|A\|^{n-1}}
$$
holds.
\end{Lem}

\subsection*{Proof of Theorem \ref{thm-4}} Let $\rho$ be the radius of the largest univalence ball of
center $0$ and contained in $f(\mathbb{B}^{n}(r)).$ Then  we have $\|f(z_{0})\|=\rho$ for some $z_{0}$ with $\|z_{0}\|=r$.
Let $[0,f(z_{0})]$ denote the segment from $0$ to $f(z_{0})$ and $\gamma$ be a curve joining $0$ and $z_{0}$ in $\mathbb{B}^{n}(r)$,
which is the preimage of $[0,f(z_{0})]$ for the mapping $f$. We use $\gamma(t)$ to denote a smooth parametrization of $\gamma$ with
$\gamma(0)=0$ and $\gamma(1)=z_{0},$ where $t\in[0,1].$


Applying \cite[Theorem 4.1 (4.2)]{PS} and Lemma \Ref{Lem-A}, we get
\beq\label{eq-y1}
\nonumber \left\|\partial_{\theta}f(z)\right\|&=&\left\|Dh(z)\theta+\overline{Dg(z)[Dh(z)]^{-1}Dh(z)\theta}\right\|\\
\nonumber &\geq&\|Dh(z)\theta\|\left(1-\big\|Dg(z)[Dh(z)]^{-1}\big\|\right)\\
\nonumber  &\geq&(1-k)\|Dh(z)\theta\|\\
\nonumber  &=&(1-k)\left\|Dh(0)\frac{[Dh(0)]^{-1}Dh(z)\theta}{\big\|[Dh(0)]^{-1}Dh(z)\theta\big\|}\right\|\big\|[Dh(0)]^{-1}Dh(z)\theta\big\|\\
\nonumber &\geq&(1-k)\frac{(1-\|z\|)^{(2n-1)\alpha+(n-3)/2}}{(1+\|z\|)^{(2n-1)\alpha-(n-3)/2}}
\min_{\xi\in\mathbb{B}^{n}}\|Dh(0)\xi\|
\eeq
which implies that
\begin{eqnarray*}
\rho&=&|f(z_{0})|=\left\|\int_{0}^{1}\frac{d}{dt}f\big(\gamma(t)\big)dt\right\|\\
&=&\int_{0}^{1}\left\|\frac{d}{dt}f\big(\gamma(t)\big)\right\|dt=\int_{0}^{1}
\left\|\partial_{\theta}f\big(\gamma(t)\big)\right\|\, |\gamma'(t)|\, dt\\
&\geq&(1-k)\min_{\theta\in\mathbb{B}^{n}}\|Dh(\gamma(0))\theta\|
\int_{0}^{1}\frac{(1-\|\gamma(t)\|)^{(2n-1)\alpha+(n-3)/2}}{(1+\|\gamma(t)\|)^{(2n-1)\alpha-(n-3)/2}}\|\,d\gamma(t)\|
\\
&\geq&(1-k)\min_{\theta\in\mathbb{B}^{n}}\|Dh(0)\theta\|\int_{0}^{r}
\frac{(1-\|z\|)^{(2n-1)\alpha+(n-3)/2}}{(1+\|z\|)^{(2n-1)\alpha-(n-3)/2}}\,d\|z\|\\
&\geq&\frac{(1-k)|\det Dh(0)|}{\|Dh(0)\|^{n-1}}\int_{0}^{r}
\frac{(1-\|z\|)^{(2n-1)\alpha+(n-3)/2}}{(1+\|z\|)^{(2n-1)\alpha-(n-3)/2}}\,d\|z\|,
\end{eqnarray*}
where $\gamma'(t)=|\gamma'(t)|\theta.$

In particular, if $n=1$, then
\begin{eqnarray*}
\rho&\geq&(1-k)\min_{\xi\in\mathbb{B}^{n}}\|Dh(0)\xi\|\int_{0}^{r}
\frac{(1-\|z\|)^{(2n-1)\alpha+(n-3)/2}}{(1+\|z\|)^{(2n-1)\alpha-(n-3)/2}}\,d\|z\|\\
&\geq&\frac{1-k}{1+k}\int_{0}^{r}\frac{(1-x)^{\alpha-1}}{(1+x)^{\alpha+1}}dx\\
&=&\frac{1-k}{2\alpha(1+k)}\left[1-\left(\frac{1-r}{1+r}\right )^{\alpha}\right].
\end{eqnarray*}
The proof of the theorem is complete.
\hfill $\Box$

\begin{lem}\label{lem-1}
Suppose that $A=(a_{ij})$ is an $n\times n$ matrix. Then
$$\left(\min_{\theta\in\partial\mathbb{B}^{n}}\|A\theta\|\right)^{n}\leq|\det A|\leq \|A\|^{n}.
$$
\end{lem}
\bpf If $A^{\ast}=(\overline{a_{ji}}),$ then the product $A^{\ast}A$
is a positive semi-definite matrix. Let $\lambda_{1},\ldots,
\lambda_{n}\,(0\leq\lambda_{1}\leq\cdots\leq\lambda_{n})$ be the $n$
eigenvalues of the matrix $A^{\ast}A$. Then
$$\sqrt{\lambda_{n}}=\max\{\|A\theta\|:\,\theta\in\partial\mathbb{B}^{n}\}~\mbox{ and }~
\sqrt{\lambda_{1}}=\min\{\|A\theta\|:\,\theta\in\partial\mathbb{B}^{n}\},
$$
which implies that
$$\|A\|^{n}\geq|\det A|=\sqrt{\Pi_{k=1}^{n}\lambda_{k}}\geq\big(\sqrt{\lambda_{1}}\big)^{n}
=\left(\min_{\theta\in\partial\mathbb{B}^{n}}\|A\theta\|\right)^{n}.
$$
The proof of the lemma is complete. \epf

\subsection*{Proof of Theorem \ref{thm-5}} In view of Lemma \ref{lem-1} and \cite[Theorem 5.1]{Pf},
$J_f$ given by \eqref{eq-ex2a} shows that
\beq\label{eq-y2} \nonumber |\det J_{f}(z)|&=&|\det
Dh(z)|^{2}\det\left(I_{n}-Dg(z)[Dh(z)]^{-1}\overline{Dg(z)[Dh(z)]^{-1}}\right)\\
\nonumber
 &\geq&|\det
Dh(z)|^{2}\min_{\theta\in\partial\mathbb{B}^{n}}\left\|\big(I_{n}-Dg(z)[Dh(z)]^{-1}\overline{Dg(z)[Dh(z)]^{-1}}\big)\theta\right\|^{n}\\
\nonumber &\geq&|\det Dh(z)|^{2}\left(1-\big\|Dg(z)[Dh(z)]^{-1}\big\|^{2}\right)^{n}\\
\nonumber &\geq&|\det Dh(z)|^{2}(1-k^{2})^{n}\\
\nonumber &=&\frac{\left|\det
\big([Dh(0)]^{-1}Dh(z)\big)\right|^{2}(1-k^{2})^{n}}{(\det[Dh(0)]^{-1})^{2}}\\
\nonumber &\geq&\frac{(1-k^{2})^{n}}{\big(\det
[Dh(0)]^{-1}\big)^{2}}\frac{\big(1-\|z\|\big)^{2n\alpha-n-1}}
{\big(1+\|z\|\big)^{2n\alpha+n+1}}.
\eeq
The proof of the theorem is complete. \hfill $\Box$

\subsection*{Proof of Theorem \ref{th-7}} By the inverse mapping theorem, we know  that $f^{-1}$ is
differentiable. Let $f^{-1}=(\sigma_{1}~\cdots~\sigma_{n})^{T}$. Then for $j,~m\in\{1,\ldots, n\}$, we use
$Df^{-1}$ and $\overline{D}f^{-1}$ to denote the two $n\times n$ matrices
$\left(\partial \sigma_{j}/\partial z_{m}\right)_{n\times n}$ and
$\left(\partial \sigma_{j}/\partial \overline{z}_{m}\right)_{n\times n}$,
respectively.

Differentiation of the equation $f^{-1}(f(z))=z $ yields the following relations
$$\begin{cases}
\displaystyle Df^{-1} Dh+\overline{D}f^{-1} Dg=I_{n},\\
\displaystyle Df^{-1} \overline{Dg}+\overline{D}f^{-1}
\overline{Dh}=0,
\end{cases}
$$
which give
\be\label{eq-2f}
\begin{cases}
\displaystyle Dh Df^{-1}=\left(I_{n}-\overline{Dg}[\overline{Dh}]^{-1}Dg[Dh]^{-1}\right)^{-1},\\
\displaystyle
Dh\overline{D}f^{-1}=-\left(I_{n}-\overline{Dg}[\overline{Dh}]^{-1}Dg[Dh]^{-1}\right)^{-1}\overline{Dg}[\overline{Dh}]^{-1}.
\end{cases}
\ee  By (\ref{eq-2f}), we get
 \beq \label{eq-1f} \nonumber \|DhDf^{-1}\|+\|Dh\overline{D}f^{-1}\| &=&
\big\|\left(I_{n}-\overline{Dg}[\overline{Dh}]^{-1}Dg[Dh]^{-1}\right)^{-1}\big\|\\
\nonumber
&&+\big\|\left(I_{n}-\overline{Dg}[\overline{Dh}]^{-1}Dg[Dh]^{-1}\right)^{-1}\overline{Dg}[\overline{Dh}]^{-1}\big\|\\
\nonumber
&\leq&\big\|\left(I_{n}-\overline{Dg}[\overline{Dh}]^{-1}Dg[Dh]^{-1}\right)^{-1}\big\|
\left(1+\|Dg[Dh]^{-1}\|\right)\\
\nonumber&\leq&\frac{1+\|Dg[Dh]^{-1}\|}
{1-\big\|\overline{Dg}[\overline{Dh}]^{-1}Dg[Dh]^{-1}\big\|}\\
&\leq&\frac{1+\|Dg[Dh]^{-1}\|}{1-\|Dg[Dh]^{-1}\|^{2}} =
\frac{1}{1-\|Dg[Dh]^{-1}\|}.
\eeq
Since $\Omega=f(\overline{\mathbb{B}^{n}(r)})$ is starlike, for each point $z_{0}\in\overline{\mathbb{B}^{n}(r)}$
and $t\in[0,1]$, we have $\varphi(t)=tf(z_{0})\in\Omega, $ where $f=(f_{1}~\cdots~f_{n})^{T}$.
Let $\gamma=f^{-1}\circ\varphi$. For any fixed $\theta\in\partial\mathbb{B}^{n}$, let
$A_{\theta}=Dg[Dh]^{-1}\theta$. By Schwarz's lemma, for
$z\in\mathbb{B}^{n}(r)$, $\|A_{\theta}(z)\|\leq\|z\|$ if
$r\in(0,1)$. The arbitrariness of $\theta\in\partial\mathbb{B}^{n}$
gives \be\label{eq-3f} \|Dg(z)[Dh(z)]^{-1}\|\leq\|z\|\leq r \ee for
$z\in\mathbb{B}^{n}(r)$. As before, by (\ref{eq-1f}) and
(\ref{eq-3f}), we obtain that
\begin{eqnarray*}
\|h(z_{0})\|
&=&\left\|\int_{0}^{1}Dh(\gamma(t))\,\frac{d}{dt}\gamma(t)\,dt\right\|\\
&=&\left\|\int_{0}^{1}Dh(\gamma(t))\left [Df^{-1}(\varphi(t))D\varphi(t)+
\overline{D}f^{-1}(\varphi(t))\overline{D\varphi(t)}\right ]\,dt\right\|\\
&\leq&\int_{0}^{1}\big(\|Dh(\gamma(t))Df^{-1}(\varphi(t))\|+\|Dh(\gamma(t))\overline{D}f^{-1}(\varphi(t))\|\big)\|D\varphi(t)\|\,dt\\
&\leq&\|f(z_{0})\|\int_{0}^{1}(1+\|Dg(\gamma(t))[Dh(\gamma(t))]^{-1}\|)\\
&&\times\left \|I_{n}-\overline{Dg(\gamma(t))}[\overline{Dh(\gamma(t))}]^{-1}
Dg(\gamma(t))[Dh(\gamma(t))]^{-1}\right \|\,dt\\
&\leq&\int_{0}^{1}\frac{1+\|Dg(\gamma(t))[Dh(\gamma(t))]^{-1}\|}
{1-\left \|\overline{Dg(\gamma(t))}[\overline{Dh(\gamma(t))}]^{-1}Dg(\gamma(t))[Dh(\gamma(t))]^{-1}\right\|}\,dt\\
&&\times\|f(z_{0})\|\\
&\leq&\|f(z_{0})\|\int_{0}^{1}\frac{1}{1-\left \|Dg(\gamma(t))[Dh(\gamma(t))]^{-1}\right \|}\,dt\\
&\leq&\frac{1}{1-r}\|f(z_{0})\|,
\end{eqnarray*}
where
$$D\varphi(t)=\left(\begin{array}{ccccc}
\ds f_{1}(z_{0})&  0  & 0 & \cdots &  0  \\[2mm]
\ds 0&   f_{2}(z_{0}) &  0&\cdots &   0\\[2mm]
 \ds \vdots & \vdots & \vdots & \cdots & \vdots \\[2mm]
 \ds  0& 0  &    \cdots &
 f_{n-1}(z_{0})&   0\\[2mm]
\ds 0&    0&   \cdots &
 0&    f_{n}(z_{0})
\end{array}\right)
$$
is a diagonal matrix.

Now we prove the second part of Theorem \ref{th-7}(a) and (b). By \cite[Theorem 5.7]{PS}, we know that
$h(\mathbb{B}^{n}(r_{0}))$ is starlike.  For $\zeta\in\mathbb{B}^{n}$, let $H(\zeta)=h(r_{0}\zeta)/r_{0}$.
Applying \cite[Theorem 2.1]{Ba} to $H$, we know that for $\zeta\in\mathbb{B}^{n}$,
$$\|H(\zeta)\|\geq\frac{\|\zeta\|}{(1+\|\zeta\|)^{2}},
$$
which implies for $z\in\mathbb{B}^{n}(r_{0})$,
\be\label{eq-y}
\|h(z)\|\geq\frac{r_{0}^{2}\|z\|}{(r_{0}+\|z\|)^{2}}.
\ee
Then Theorem \ref{th-7} (a) follows from (\ref{eq-y}), and  Theorem \ref{th-7} (b) easily
follows from Theorem \ref{th-7}(a). The proof of the theorem is complete.
\hfill $\Box$


\subsection*{Proof of Theorem \ref{thm-6}} We first prove the sufficiency of part (a). Without loss of generality, we
assume that
\be\label{eq-1.9}
\|Dh(z)\|\leq K|\det Dh(z) |^{\frac{1}{n}} ~\mbox{ for $z\in\mathbb{B}^{n}$},
\ee
where $K\geq1$ is a constant.

As in the proof of Theorem \ref{thm-5}, (\ref{eq-1.9}) and Lemma \ref{lem-1}, for $z\in\mathbb{B}^{n}$,
we have
\begin{eqnarray*}
|\det J_{f}(z)|
&\geq&|\det Dh(z)|^{2}(1-c^{2})^{n}
\end{eqnarray*}
so that
$$ |\det Dh(z)|^{\frac{1}{n}} \leq \frac{|\det J_{f}(z) |^{\frac{1}{2n}}}{\sqrt{1-c^2}}.
$$
Moreover,
$$\Lambda_{f}(z)=\max_{\theta\in\partial\mathbb{B}_{\mathbb{R}}^{2n}}\|J_{f}(z)\theta\|\leq\|Dh(z)\|
\left(1+\big\|Dg(z)[Dh(z)]^{-1}\big\|\right)\leq\|Dh(z)\|(1+c),
$$
which by the last inequality gives that
\be\label{eq-u1}
\Lambda_{f}(z)\leq K\sqrt{\frac{1+c}{1-c}}\, |\det J_{f}(z)|^{\frac{1}{2n}}
\ee
and hence, $f$ is a quasiregular mapping.

Next we  prove the necessity of part (a). We assume that for $z\in\mathbb{B}^{n}$,
\be\label{eq-1.91}
\Lambda_{f}(z)\leq K_{1}|\det J_{f}(z) |^{\frac{1}{2n}},
\ee
where $K_{1}\geq1$ is a constant.

As in the proof of Theorem \ref{thm-5}, for $z\in\mathbb{B}^{n}$, by calculations and Lemma \ref{lem-1}, we get
\begin{eqnarray*}
|\det J_{f}(z)|&=&|\det
Dh(z)|^{2}\left|\det\left(I_{n}-Dg(z)[Dh(z)]^{-1}\overline{Dg(z)[Dh(z)]^{-1}}\right )\right|\\
&\leq&|\det
Dh(z)|^{2}\left\|I_{n}-Dg(z)[Dh(z)]^{-1}\overline{Dg(z)[Dh(z)]^{-1}}\right\|^{n}\\
&\leq&|\det Dh(z)|^{2}\left(1+c^{2}\right)^{n}
\end{eqnarray*}
so that
$$ |\det Dh(z)|^{\frac{1}{n}} \geq \frac{|\det J_{f}(z) |^{\frac{1}{2n}}}{\sqrt{1+c^2}}.
$$
Furthermore,
 $$\Lambda_{f}(z)=\max_{\theta\in\partial\mathbb{B}_{\mathbb{R}}^{2n}}\|J_{f}(z)\theta\|\geq\|Dh(z)\|
\left(1-\big\|Dg(z)[Dh(z)]^{-1}\big\|\right)\geq\|Dh(z)\|(1-c),
$$
which, by \eqref{eq-1.91}, implies that
 $$\|Dh(z)\|(1-c)\leq\Lambda_{f}(z)\leq
K_{1}|\det J_{f}(z) |^{\frac{1}{2n}}\leq K_{1}\sqrt{1+c^{2}}\, |\det
Dh(z) |^{\frac{1}{n}}.
$$
Hence
$$\|Dh(z)\|\leq\frac{K_{1}\sqrt{1+c^{2}}}{1-c}|\det Dh(z) |^{\frac{1}{n}},
$$
which shows that $h$ is a quasiregular mapping.

Now we prove part (b). By (\ref{eq-u1}), we know that
$f$ is a pluriharmonic $K_{2}$-quasiregular mapping, where
$K_{2}=K\sqrt{\frac{1+c}{1-c}}$. Applying \cite[Theorem 6]{HG1}, we
know that $f(\mathbb{B}^{n})$ contains a univalent ball with the
radius $R$ with
$$R\geq\frac{k_{n}\pi}{8m}\left(\frac{k_{n}\pi}{4K_{2}\log(1/(1-k_{n}))}\right)^{4n-1},
$$
where $m\approx4.2$ and $0<k_{n}<1$ is a unique root such that
$$4n\log\frac{1}{1-k_{n}}=(4n-1)\frac{k_{n}}{1-k_{n}}.
$$
The proof of the theorem is complete. \hfill $\Box$

\subsection*{Acknowledgements}
The research of the second author was supported by the project RUS/RFBR/P-163 under
Department of Science \& Technology and Russian Foundation for Basic Research and this
author is currently on leave from the Department of Mathematics,
Indian Institute of Technology Madras, Chennai-600 036, India.

\normalsize

\end{document}